\documentclass[11pt]{article}
\usepackage{amsfonts,amsmath,amssymb}
\usepackage{hyperref}
\usepackage{graphicx}
\usepackage{color}

\newcommand{\halmos}{\rule{1ex}{1.4ex}}
\def \qed {\nopagebreak{\hspace*{\fill}$\halmos$\medskip}}

\newtheorem{theorem}{Theorem}[section]
\newtheorem{proposition}[theorem]{Proposition}
\newtheorem{corollary}[theorem]{Corollary}
\newtheorem{conjecture}[theorem]{Conjecture}
\newtheorem{lemma}[theorem]{Lemma}

\newtheorem{remark}[theorem]{Remark}
\newtheorem{defi}[theorem]{Definition}

\newcommand{\bt}{\begin{theorem}}
\newcommand{\et}{\end{theorem}}
\newcommand{\bl}{\begin{lemma}}
\newcommand{\el}{\end{lemma}}
\newcommand{\bprop}{\begin{proposition}}
\newcommand{\eprop}{\end{proposition}}
\newcommand{\bcor}{\begin{corollary}}
\newcommand{\ecor}{\end{corollary}}
\newcommand{\br}{\begin{remark}\rm}
\newcommand{\er}{\end{remark}}
\newcommand{\bcon}{\begin{conjecture}}
\newcommand{\econ}{\end{conjecture}}
\newcommand{\bd}{\begin{defi}}
\newcommand{\ed}{\end{defi}}

\newcommand{\be}{\begin{equation}}
\newcommand{\ee}{\end{equation}}

\newcommand{\Di}{{\cal D}}

\newcommand{\Hi}{{\cal H}}

\newcommand{\Ki}{{\cal K}}

\newcommand{\Ri}{{\cal R}}

\newcommand{\Wi}{{\cal W}}

\newcommand{\eps}{\epsilon}

\newcommand{\R}{{\mathbb R}}
\newcommand{\N}{{\mathbb N}}
\newcommand{\Z}{{\mathbb Z}}

\renewcommand{\P}{{\mathbb P}}
\newcommand{\E}{{\mathbb E}}

\newcommand{\Asto}[1]{\underset{{#1}\to\infty}{\Longrightarrow}}

\newcommand{\Rc}{R^2_{\rm c}}

\begin{document}

\makeatletter\@addtoreset{equation}{section}
\makeatother\def\theequation{\thesection.\arabic{equation}}

\renewcommand{\labelenumi}{{(\roman{enumi})}}

\title{Brownian Web and Oriented Percolation: Density Bounds}
\author{Anish Sarkar$^{\,1}$ \and Rongfeng Sun$^{\,2}$}

\maketitle

\footnotetext[1]{Theoretical Statistics and Mathematics Unit, Indian Statistical Institute, New Delhi,
7 S.~J.~S.\ Sansanwal Marg, New Delhi 110016, India. Email: anish@isid.ac.in}

\footnotetext[2]{Department of Mathematics, National University of Singapore, 10 Lower Kent Ridge Road, 119076 Singapore. Email:
matsr@nus.edu.sg}

\begin{abstract}

In a recent work~\cite{SS11}, we proved that under diffusive scaling, the collection of rightmost infinite open paths in a supercritical oriented percolation
configuration on the space-time lattice $\Z^2$ converges in distribution to the Brownian web. In that proof, the FKG inequality played an important role in establishing a density bound, which is a part of the convergence criterion for the Brownian web formulated by Fontes et al in~\cite{FINR04}. In this note, we illustrate how an alternative convergence criterion formulated by Newman et al in~\cite{NRS05} can be verified in this case, which involves a dual density bound that can be established without using the FKG inequality. This alternative approach is in some sense more robust. We will also show that the spatial density of the collection of rightmost infinite open paths starting at time $0$ decays in time as $\frac{2+o(1)}{\sigma\sqrt{\pi t}}$ for some $\sigma>0$.

\end{abstract}

\vspace{6pt}

\section{Introduction}

We first briefly recall the basic setup from~\cite{SS11}. Since this note is meant to complement~\cite{SS11}, we will refer the reader to~\cite{SS11} for
many details.

Let $\Z^2_{\rm even}:=\{(x,i)\in \Z^2: x+i \mbox{ is even}\}$ be a space-time lattice, with oriented edges leading from $(x,i)$ to $(x\pm 1, i+1)$ for all $(x,i)\in \Z^2_{\rm even}$. For a fixed parameter $p\in [0,1]$, independently each oriented edge is open with probability $p$, and closed with probability $1-p$. The random
configuration of open and closed oriented edges defines the oriented percolation model on $\Z^2_{\rm even}$. We will use $\P_p$ and $\E_p$ to denote respectively probability and expectation for this product probability measure on edge configurations with parameter $p$, with $p$ omitted when there is no confusion.

For each $z\in\Z^2_{\rm even}$, we will let $C_z$ denote the open cluster at $z$, which contains all sites in $\Z^2_{\rm even}$ that can be reached from $z$ by following open oriented edges. We will denote by $|C_z|$ the cardinality of $C_z$. When $|C_z|=\infty$, $z$ is called a {\em percolation point}, and we denote the set of percolation points by $\Ki$. In the supercritical regime $p\in (p_c, 1]$, $\Ki$ is almost surely an infinite set.

For each $z=(x,i)\in \Ki$, there is a well-defined {\em rightmost infinite open path} starting from $z$, which we denote by $\gamma_z$. More precisely, $\gamma_z$ can be taken as a mapping from $\{i, i+1, \cdots\}$ to $\Z$ such that $\gamma_z(i)=x$, there is an oriented edge from $(\gamma_z(j), j)$ to $(\gamma_z(j+1), j+1)$ which is open for all $j\geq i$, and if $\pi$ is any other infinite open path starting from $z$, then $\gamma_z(j)\geq \pi(j)$ for all $j\geq i$. We will identify $\gamma_z$ with the continuous path (i.e., function) which is defined on $[i, \infty)$ by linearly interpolating between $\gamma_z$ at consecutive integer times. If $z=(x,i)\in\Z^2_{\rm even}$ is not a percolation point, then we define
\be\label{gammaz}
\gamma_z:=\gamma_{z'}, \quad \mbox{where } z'=(y,i) \mbox{ with } y=\max\{u\leq x : (u,i) \in \Ki\},
\ee
which is the rightmost infinite open path starting from $(-\infty, x]\times\{i\}$. We are interested in the collection of all rightmost infinite open paths in the supercritical oriented percolation configuration (for a fixed $p\in (p_c, 1)$ which we assume from now on):
\be\label{Gamma}
\Gamma := \{\gamma_z: z\in\Z^2_{\rm even}\}= \{\gamma_z : z\in \Ki\}.
\ee
As in~\cite{SS11}, we denote the space of {\em continuous paths} by $\Pi$, equipped with a suitable metric $d$, and we denote the space of {\em compact subsets of $\Pi$} by $\Hi$, equipped with the Hausdorff metric $d_\Hi$ induced by $d$. It can be shown that $\overline{\Gamma}$, the closure of $\Gamma$ in $(\Pi, d)$,
is almost surely compact, and $\overline{\Gamma}\backslash\Gamma$ contains only trivial paths which arise from the compactification of $\R^2$ and $\Pi$.
See~\cite{SS11} for more details. Therefore we will regard $\overline{\Gamma}$ as an $(\Hi, d_\Hi)$-valued random variable.  To simply notation at the cost of a slight abuse,  in what follows, we will write $\Gamma$ instead of $\overline{\Gamma}$, and we will omit taking closure of sets of points or paths explicitly when such closure only adds trivial elements.

Let $o$ denote the origin. It is known from results of Durret~\cite{D84} and Kuczek~\cite{K89} that there exists $\alpha=\alpha(p)>0$ and $\sigma=\sigma(p)>0$ such that
$\frac{\gamma_o(n)-\alpha n}{\sigma \sqrt{n}}$ converges in distribution to a standard normal as $n\to\infty$. For $\eps>0$, let $S_\eps:\R^2\to\R^2$ denote the
shearing and diffusive scaling map with
\be\label{Smap}
S_{\eps} (x,t) :=\Big(\frac{\sqrt{\eps}}{\sigma}(x-\alpha t), \eps t\Big).
\ee
When $t$ is understood to be a time, we will define $S_{\eps}t:=\eps t$. By identifying a continuous path $\pi: [\sigma_\pi,\infty)\to\R$ (with starting time denoted by $\sigma_\pi$) with its graph in $\R^2$, we can define $S_\eps \pi$ by applying $S_\eps$ to each point on the graph of $\pi$. Similarly if $K\subset \Pi$, we can define $S_\eps K$ by applying $S_\eps$ to each element of $K$.

In~\cite{SS11}, we proved the following convergence result, which verifies a conjecture of Wu and Zhang~\cite{WZ08}.
\bt\label{T:main} Let $p\in (p_c, 1)$, and let $\sigma,\alpha$ be as in (\ref{Smap}). As $\eps\downarrow 0$, $S_{\eps}\Gamma$ converges in distribution to the standard Brownian web $\Wi$ as $(\Hi, d_\Hi)$-valued random variables.
\et
Loosely speaking, the Brownian web $\Wi$ is the set of coalescing Brownian motions starting from every point in the space-time plane $\R^2$, so that Theorem~\ref{T:main} effectively asserts that after centering and diffusively rescaling, the paths in $\Gamma$ converge in distribution to a collection of coalescing Brownian motions. The rigorous formulation is more subtle than this heuristic description, and we refer the reader to~\cite{SS11} for details.

Theorem~\ref{T:main} was proved in~\cite{SS11} by verifying the following convergence criterion for the Brownian web, which was proposed by Fontes et al in~\cite{FINR04}.

Let $K\in \Hi$. For $\pi \in K$, recall that $\sigma_\pi$ denotes the starting time of $\pi$. For $t>0$ and $t_0,a,b\in\R$ with $a<b$, let
\be\label{eta}
\eta_K(t_0,t;a,b) := \big|\{\pi(t_0+t)\, :\, \pi \in K \mbox{ with } \sigma_\pi\leq t_0 \mbox{ and } \pi(t_0)\in [a,b]  \}\big|,
\ee
which counts the number of distinct points on $\R\times\{t_0+t\}$ touched by some path in $K$ which also touches $[a,b]\times\{t_0\}$.

A random variable $\cal X$ taking values in $\Hi$ is said to have {\em non-crossing paths} if a.s.\ there exist no $\pi, \tilde\pi \in {\cal X}$ such that
$$
(\pi(t)-\tilde\pi(t))(\pi(s)-\tilde\pi(s))<0 \qquad \mbox{for some } s,t\geq \sigma_\pi \vee \sigma_{\tilde\pi}.
$$
In our case, $\Gamma$, and hence also $S_\eps\Gamma$, have non-crossing path.

\bt\label{T:convct1}{\bf \cite{FINR04}} Let $({\cal X}_n)_{n\in\N}$ be a sequence of $\Hi$-valued random variables with
non-crossing paths. If the following conditions are satisfied, then ${\cal X}_n$ converges in distribution to the standard Brownian web $\Wi$.
\begin{itemize}
\item[{\rm \,(I)\ }] Let $\Di$ be a deterministic countable dense subset of $\R^2$. Then there exist $\pi^y_n \in {\cal X}_n$ for $y\in\Di$ such that, for each
finite collection $y_1, y_2, \ldots, y_k\in \Di$, $(\pi^{y_1}_n, \ldots, \pi^{y_k}_n)$ converge in distribution as $n\to\infty$ to a collection of coalescing Brownian motions starting at $(y_1,\ldots, y_k)$.

\item[{\rm (B1)}] For all $t>0$, $\limsup_{n\to\infty} \sup_{(a,t_0)\in\R^2} \P(\eta_{{\cal X}_n}(t_0, t; a, a+\eps)\geq 2) \to 0$ as $\eps\downarrow 0$.
\item[{\rm (B2)}] For all $t>0$, $\eps^{-1}\limsup_{n\to\infty} \sup_{(a,t_0)\in\R^2} \P(\eta_{{\cal X}_n}(t_0, t; a, a+\eps)\geq 3) \to 0$ as $\eps\downarrow 0$.
\end{itemize}
\et
As shown in \cite[Prop.~B.2]{FINR04}, condition (I) and the non-crossing property imply that $({\cal X}_n)_{n\in\N}$ is a tight sequence of $\Hi$-valued
random variables. Condition (I) also guarantees that any subsequential weak limit of $({\cal X}_n)_{n\in\N}$ contains as many paths as, possibly more than,
the Brownian web $\Wi$. Conditions (B1) and (B2) are density bounds which rule out the presence of extra paths other than the Brownian web paths in any
subsequential weak limit.

In~\cite{SS11}, the conditions in Theorem~\ref{T:convct1} were verified for the sequence of $\Hi$-valued random variables $S_\eps \Gamma$, with $\eps\downarrow 0$. The key idea was the approximation of each path in $\Gamma$ by a percolation exploration cluster. Each exploration cluster evolves in a Markovian way, and different exploration clusters evolve independently before they intersect. Condition (I) was verified by establishing an invariance principle for each exploration cluster and then showing that distinct exploration clusters coalesce when they intersect. Condition (B1) then follows easily from (I). Using the approximation by exploration clusters, condition (B2) was verified by applying the FKG inequality for the oriented percolation edge configuration. However, the fact that one could apply the FKG inequality to prove (B2) was not obvious at first, and it was remarked at the end of Section~1.4 in~\cite{SS11} that there is an alternative approach to proving the convergence of $S_\eps\Gamma$ to the Brownian web, which does not rely on the FKG inequality and is in a sense more robust. The goal of this note is to provide the details for this alternative approach.

We will use an alternative convergence criterion developed in~\cite{NRS05}, which replaces conditions (B1) and (B2) by a density bound on a dual counting variable
\be\label{etahat}
\hat \eta_K(t_0,t;a,b) := \big|(a,b) \cap \{\pi(t_0+t)\, :\, \pi \in K \mbox{ with } \sigma_\pi\leq t_0 \}\big|,
\ee
which counts the number of distinct points on $(a,b)\times\{t_0+t\}$ touched by some path in $K$ starting before or at time $t_0$.
Here is the convergence criterion formulated in~\cite[Theorem 1.4]{NRS05}, where we have removed the tightness condition therein, which is redundant for $\Hi$-valued random variables with non-crossing paths.

\bt\label{T:convct2}{\bf \cite{NRS05}} Let $({\cal X}_n)_{n\in\N}$ be a sequence of $\Hi$-valued random variables with
non-crossing paths. If condition {\rm (I)} from Theorem~\ref{T:convct1} and the following condition hold, then ${\cal X}_n$ converges in distribution to the standard Brownian web $\Wi$:
\begin{itemize}
\item[{\rm (E)}] If ${\cal X}$ is any subsequential weak limit of $({\cal X}_n)_{n\in\N}$, then for all $t>0$ and $t_0,a,b\in\R$ with $a<b$,
we have $\E[\hat\eta_{\cal X}(t_0,t;a,b)] \leq \E[\hat\eta_{\Wi}(t_0,t;a,b)] = \frac{b-a}{\sqrt{\pi t}}$.
\end{itemize}
\et
{\bf Remark.} We point out that condition $(B'_1)$ in~\cite[Theorem 1.4]{NRS05}, which is a variant of (B1), is in fact redundant in the formulation of that theorem. This is because (B1) and (B2) together ensure that any subsequential weak limit contains no more paths than in the Brownian web, which is now guaranteed by condition (E).
For further details, see~\cite[Theorem 4.2]{FINR04} and the remark afterwards, and the discussion before \cite[Theorem 1.4]{NRS05}. Condition $(B'_1)$ was however used in \cite{NRS05} to verify condition (E), because paths there can cross.
\medskip

To verify condition (E), we will follow a strategy developed in \cite{NRS05}. First we show that the density at time $t$ of the set of coalescing paths in $\Gamma$, starting before or at time $0$, decays with the order of $O(1/\sqrt{t})$, which guarantees that at any time $\delta>0$ on the diffusive space-time scale, only a locally finite number of paths remain. This bound is conceptually easier than (B2) because it is equivalent to a bound on the expectation of the random variable $\eta_{{\cal X}_n}(t_0, t; a, a+\delta)$ in (B2). We can then use the Brownian approximation given by condition (I) to refine the density upper bound to establish (E). In the last step, we will need to establish some asymptotic Markovian properties for paths in $\Gamma$, which is enjoyed by the Brownian web. We remark that apart from~\cite{NRS05}, another instance where the convergence criterion in Theorem~\ref{T:convct2} was used to prove convergence to the Brownian web is a generalized drainage network model studied recently in~\cite{CV11}.

Using Theorem~\ref{T:main}, we can in turn refine the density bound for paths in $\Gamma$, which was used in the verification of condition (E).
\bt\label{T:density} Let $p\in (p_c, 1)$, and let $\sigma$ be as in (\ref{Smap}). Let $\Gamma_0(n):=\{\gamma_{(2x,0)}(n) : x\in\Z\}$, which is supported on the subset
of $\Z$ with the same parity as $n$. Then as $n\uparrow\infty$,
\be\label{density}
\P(\Gamma_0(n)\cap \{0,1\}\neq \emptyset) = \frac{2+o(1)}{\sigma \sqrt{\pi n}}.
\ee
\et

We will establish condition (E) and Theorem~\ref{T:density} in Sections~\ref{S:E} and \ref{S:D}, respectively.

\section{Verification of Condition (E)}\label{S:E}
We will verify condition (E) in Theorem~\ref{T:convct2} for $S_\eps\Gamma$, with $\eps\downarrow 0$, by following the strategy developed in~\cite{NRS05}.
Recall that condition (I) was verified in \cite[Prop.~3.3]{SS11}.

In what follows, for any $s<t$ and $K\in \Hi$, we will denote
\be\label{truncation}
\begin{aligned}
K_s& :=\{\pi\in K: \sigma_\pi =s\}, \qquad & K_{s^-} &:= \{\pi \in K : \sigma_\pi \leq s\},  \\
K(s) & := \{\pi(s) : \pi \in K\},  \qquad &  K^t& := \{\pi^t : \pi \in K\},
\end{aligned}
\ee
where $\pi^t$ denotes the path obtained from $\pi$ by restricting $\pi$ to the time interval $[t, \infty)$. We will also denote $K_{s^-}^t:=(K_{s^-})^t$.

First we recall from~\cite[Lemma 6.1]{NRS05} that condition (E) can be replaced by condition
\begin{itemize}
\item[{\rm (E')}] For any $t_0\in\R$, if ${\cal Z}$ is a subsequential weak limit of $({\cal X}_{n, t_0^-})_{n\in\N}$, where
${\cal X}_{n, t_0^-}:=({\cal X}_n)_{t_0^-}$, then for all $t>0$ and $a<b$, we have
\be\label{Zbound}
\E[\hat\eta_{{\cal Z}}(t_0,t;a,b)] \leq \E[\hat\eta_{\Wi}(t_0,t;a,b)] = \frac{b-a}{\sqrt{\pi t}}.
\ee
\end{itemize}
Condition (E') simplifies (E) by effectively singling out the subset of paths in ${\cal X}$ starting before or at time $t_0$, which are the only relevant paths for
verifying condition (E).

Returning to our setting, given $t_0\in\R$, let ${\cal X}_n:= S_{\eps_n}\Gamma$ for a sequence $\eps_n\downarrow 0$ such that
${\cal X}_{n, t_0^-}$ converges weakly to ${\cal Z}$. To verify (E'), it suffices to prove (\ref{Zbound}), which will follow from the next two lemmas.

\bl\label{L:Zlocfin} Let ${\cal Z}$ be the weak limit of ${\cal X}_{n, t_0^-}$, with ${\cal X}_n := S_{\eps_n}\Gamma$, for a sequence $\eps_n\downarrow 0$. Then for any $\delta>0$, ${\cal Z}(t_0+\delta)$ is a.s.\ a locally finite subset of $\R$.
\el

\bl\label{L:ZBrownian} Let $\cal Z$ be as in Lemma~\ref{L:Zlocfin}. Then for any $\delta>0$, ${\cal Z}^{t_0+\delta}$ has the same distribution as that of
$\Wi_{t_0+\delta, \cal Z}:=\{\pi \in \Wi: \sigma_\pi = t_0+\delta, \pi(t_0+\delta)\in {\cal Z}(t_0+\delta)\}$, where $\Wi$ is a standard Brownian web independent of ${\cal Z}$.
\el
Lemma~\ref{L:ZBrownian} implies that
$$
\E[\hat\eta_{{\cal Z}}(t_0,t;a,b)] = \E\big[\hat\eta_{{\cal Z}^{t_0+\delta}}(t_0+\delta,t-\delta;a,b)\big] \leq \E[\hat\eta_{{\cal W}}(t_0+\delta,t-\delta;a,b)]
=\frac{b-a}{\sqrt{\pi (t-\delta)}},
$$
from which (\ref{Zbound}) follows by letting $\delta\downarrow 0$.
\medskip

Before we proceed to prove Lemmas~\ref{L:Zlocfin}--\ref{L:ZBrownian}, we first briefly recall the exploration cluster and some of its basic properties developed in~\cite{SS11}, which we will also need to use here. Let $z=(x,t) \in\Z^2_{\rm even}$. The percolation exploration cluster $C_z(n)$ at time $n\geq t$ is the minimal set of edges one need to examine to find the rightmost open path $l_z^n:=(l_z^n(i))_{t\leq i\leq n}$ connecting $(-\infty, x]\times\{t\}$ to $\Z\times \{n\}$. The open path $l_z^n$ forms the left boundary of the exploration cluster $C_z(n)$. The right boundary of $C_z(n)$ is defined by
$$
r^n_z(j) := \max\{y\in\Z : (-\infty, x]\times\{t\} \to (y,j)\}, \qquad j=t, t+1, \ldots, n,
$$
which is the path of the rightmost vertex that can be reached at each time $j$ by an open oriented path starting from $(-\infty, x]$ at time $t$. Note that the definition
of $r^n_z(j)$ for each $j\geq t$ is independent of $n\geq j$, and hence we will just work with $r_z:=(r_z(i))_{i\geq t}$, with $r_z(j)=r^n_z(j)$ for $n\geq j\geq t$.
The boundaries of the exploration cluster, $(l^n_z, r_z)$, serve as approximations of $\gamma_z$, and it was shown in~\cite[Prop.~2.2]{SS11} that $(S_\eps \gamma_z, S_\eps r_z)$ converge in distribution to $(B, B)$ for a standard Brownian motion $B$ as $\eps\downarrow 0$. Another important fact established in~\cite[Lemma~3.1]{SS11} is
that, for $z_1=(x_1,0)\neq z_2=(x_2,0)\in\Z^2_{\rm even}$, $r_{z_1}=(r_{z_1}(n))_{n\geq 0}$ and $r_{z_2}=(r_{z_2}(n))_{n\geq 0}$ coalesce at the first time $\tau$ when
$r_{z_1}(\tau)=r_{z_2}(\tau)$, $\gamma_{z_1}$ and $\gamma_{z_2}$ must coalesce before or at time $\tau$, and the exploration clusters $C_{z_1}(n)$ and $C_{z_2}(n)$ are disjoint for $n<\tau$.

\subsection{Proof of Lemma~\ref{L:Zlocfin}}
We will prove Lemma~\ref{L:Zlocfin} via the following bound on the rate of coalescence for paths in $\Gamma$ starting at time $0$, as well as for the right
boundaries of the associated exploration clusters. Recall the notation introduced in (\ref{truncation}).

\bl\label{L:densitydecay} Let ${\cal R}:=\{r_z: z\in\Z^2_{\rm even}\}$ denote the set of right boundaries of all exploration clusters. Then
there exists $C>0$ such that for all $n\in\N$,
\be\label{ddecay}
\P(\Gamma_0(n)\cap \{0,1\}\neq\emptyset) \leq \P({\cal R}_0(n)\cap \{0,1\}\neq\emptyset) \leq \frac{C}{\sqrt n}.
\ee
\el
{\bf Proof.} Without loss of generality, we may assume that $n$ is even so that $\{0,1\}$ in (\ref{ddecay}) can be replaced by $\{0\}$.
Since paths in $\Gamma_0$ are non-crossing and coalesce when they intersect, by translation invariance, we have
$$
\begin{aligned}
\P(0\in \Gamma_0(n)) & = \sum_{x\in\Z} \P\big(\gamma_{(2x,0)}(n)=0, \gamma_{(2x+2,0)}(n)\geq 2\big) \\
& = \sum_{x\in\Z} \P\big(\gamma_{(0,0)}(n)=-2x, \gamma_{(2,0)}(n)\geq -2x+2\big) = \P(\gamma_{(0,0)}(n)<\gamma_{(2,0)}(n)),
\end{aligned}
$$
which is the probability that $\gamma_{(0,0)}$ and $\gamma_{(2,0)}$ do not coalesce by time $n$. By \cite[Lemma 3.1]{SS11} and the remark following it, paths in ${\cal R}_0$ are also non-crossing and coalesce when they intersect. Therefore the same argument implies that
$$
\P(0\in {\cal R}_0(n)) = \P(r_{(0,0)}(n) < r_{(2,0)}(n)).
$$
Also by \cite[Lemma~3.1]{SS11}, $\gamma_{(0,0)}$ and $\gamma_{(2,0)}$ must coalesce before $r_{(0,0)}$ and $r_{(2,0)}$ coalesce. Therefore
the first inequality in (\ref{ddecay}) holds.

To prove the second inequality in (\ref{ddecay}), by translation invariance, it suffices to show that
\be\label{limsupR0}
\limsup_{n\to\infty} \E\big[|{\cal R}_0(n)\cap [0, 2L-1]|\big] <\infty, \qquad \mbox{where } L=\lceil \sqrt{n}\rceil.
\ee
For $k\in \Z$, let ${\cal R}_0^k :=\{r_{(2x,0)}: x\in [kL, (k+1)L-1]\cap\Z\}$. Then by translation invariance,
$$
\begin{aligned}
\E\big[|{\cal R}_0(n)\cap [0, 2L-1]|\big] & \leq \sum_{k\in\Z} \E\big[|{\cal R}^k_0(n)\cap [0, 2L-1]|\big] \\
& = \sum_{k\in\Z} \E\big[|{\cal R}^0_0(n)\cap [-2kL, -2kL+2L-1]|\big] = \E\big[|{\cal R}^0_0(n)|\big].
\end{aligned}
$$
Therefore it suffices to show that $\E[|{\cal R}^0_0(n)|]$ is uniformly bounded in $n\in\N$. We can write
\be\label{R00}
\E\big[|{\cal R}^0_0(n)|\big] = \sum_{k=1}^\infty \P\big(|{\cal R}^0_0(n)|\geq k \big) \leq 1 + 2\sum_{k=1}^\infty \P\big(|{\cal R}^0_0(n)|\geq 2k \big).
\ee
On the event $|{\cal R}^0_0(n)|\geq 2k$ for some $k\in\N$, there exist $0\leq 2x_1 <2x_2<\cdots <2x_{2k} \leq 2L-2$ such that
$r_{(2x_1, 0)}, \cdots, r_{(2x_{2k},0)}$ are mutually disjoint on the time interval $[0, n]$. Again by \cite[Lemma~3.1]{SS11}, the realization
of $(r_{(2x_i,0)})_{1\leq i\leq 2k}$ up to time $n$ is determined by $2k$ non-intersecting exploration clusters,
and hence the realization of the $k$ pairs
$(r_{(2x_{2i-1},0)}, r_{(2x_{2i},0)})_{1\leq i\leq k}$ up to time $n$ is determined by $k$ disjoint sets of edges. Therefore if
$$
D_n:=\{\exists\ 0\leq 2y_1<2y_2\leq 2L-2 \mbox{ such that } r_{(2y_1, 0)} \mbox{ and } r_{(2y_2, 0)} \mbox{ do not coalesce by time } n\},
$$
then on the event $|{\cal R}^0_0(n)|\geq 2k$, $D_n$ occurs disjointly $k$ times. Therefore by Reimer's inequality~\cite{R00} for disjoint occurrences of events on a product probability space of Bernoulli random variables, we have
$$
\P(|{\cal R}^0_0(n)|\geq 2k) \leq \P(D_n)^k.
$$
Note that $D_n$ equals the event that $r_{(0,0)}$ and $r_{(2L-2,0)}$ do not coalesce by time $n$. Since $S_{n^{-1}}(2L-2,0)\to (2\sigma^{-1}, 0)$ as $n\to\infty$, by
the invariance principle for a pair of exploration clusters established in \cite[Prop.~3.2]{SS11}, the time of coalescence of $r_{(0,0)}$ and $r_{(2L-2,0)}$ divided by $n$ converges in distribution to the time of intersection of two independent Brownian motions starting respectively at $(0,0)$ and $(2\sigma^{-1},0)$. In particular, $\lim_{n\to\infty}\P(D_n)<1$.
Therefore the RHS of (\ref{R00}) is uniformly bounded in $n\in\N$. This proves the second inequality in (\ref{ddecay}).
\qed
\bigskip

\noindent
{\bf Proof of Lemma~\ref{L:Zlocfin}.} By Skorohod's representation theorem~\cite[Theorem 6.7]{B99}, we can couple ${\cal X}_{n, t_0^-}$, $n\in\N$, and $\cal Z$ on the same probability space such that ${\cal X}_{n, t_0^-}\to {\cal Z}$ a.s., where the convergence is in the space of compact sets of paths $(\Hi, d_\Hi)$. This implies that for any $\delta>0$, ${\cal X}_{n, t_0^-}(t_0+\delta)\to {\cal Z}(t_0+\delta)$ a.s.\ in the space of compact subsets of $[-\infty, \infty]$, the compactification of $(-\infty, \infty)$, equipped with the Hausdorff topology. In particular, for any $a<b$, a.s.\
$$
|{\cal Z}(t_0+\delta) \cap (a,b)| \leq \liminf_{n\to\infty} |{\cal X}_{n, t_0^-}(t_0+\delta)\cap (a,b)|,
$$
and hence by Fatou's lemma,
\be\label{Zt0ab}
\E[|{\cal Z}(t_0+\delta) \cap (a,b)|] \leq \liminf_{n\to\infty} \E[|{\cal X}_{n, t_0^-}(t_0+\delta)\cap (a,b)|].
\ee
Recall that ${\cal X}_{n}=S_{\eps_n}\Gamma$, and denote $\beta_n:=\eps_n^{-1}(t_0+\delta)$. Then
$$
\begin{aligned}
\E[|{\cal X}_{n, t_0^-}(t_0+\delta)\cap (a,b)|] & = \E\big[\big|\Gamma_{(\eps_n^{-1}t_0)^-}(\beta_n) \cap \big(\alpha\beta_n+a\sigma\eps_n^{-1/2},
\alpha\beta_n+b\sigma\eps_n^{-1/2}\big)\big|\big] \\
& \leq \E\big[\big|\Gamma_{\lfloor \eps_n^{-1}t_0\rfloor}(\lfloor \beta_n\rfloor) \cap \big(\alpha\beta_n+a\sigma\eps_n^{-1/2}-1,
\alpha\beta_n+b\sigma\eps_n^{-1/2}+1 \big)\big|\big] \\
& \leq \E\big[\big|\Gamma_{0}(\lfloor \beta_n\rfloor-\lfloor \eps_n^{-1}t_0\rfloor) \cap \big(\alpha\beta_n+a\sigma\eps_n^{-1/2}-1,
\alpha\beta_n+b\sigma\eps_n^{-1/2}+1 \big)\big|\big] \\
& \leq \big(2+\frac{\sigma}{2}(b-a)\eps_n^{-1/2}\big)\ \P\big(\Gamma_{0}(\lfloor \beta_n\rfloor-\lfloor \eps_n^{-1}t_0\rfloor)\cap \{0,1\}\neq\emptyset\big) \\
& \leq C \frac{\big(2+\frac{\sigma}{2}(b-a)\eps_n^{-1/2}\big)}{\sqrt{\lfloor \eps_n^{-1}(t_0+\delta)\rfloor-\lfloor \eps_n^{-1}t_0\rfloor}},
\end{aligned}
$$
which has a bounded limit as $n\to\infty$. In the first inequality above, we used the fact that $\Gamma_{s^-}^t=\Gamma_{\lfloor s\rfloor}^t$ for any
$s\in\R$ and $t>s$ because paths in $\Gamma$ coalesce when they intersect. In the next two inequalities, we used the translation invariance
of $\Gamma$ under shifts by vectors in $\Z^2_{\rm even}$, while in the last inequality we used Lemma~\ref{L:densitydecay}. This proves that the RHS of (\ref{Zt0ab}) is finite for all $a<b$, and hence ${\cal Z}(t_0+\delta)$ is a.s.\ locally finite.
\qed
\bigskip

\noindent
{\bf Remark.} Note that when bounding the RHS of (\ref{R00}), we applied Reimer's inequality, which is as strong and delicate as the FKG inequality. However, even if Reimer's inequality was not available, we believe it would still be much easier to show that the RHS of (\ref{R00}) is bounded (which is what we need here) than to show $\P\big(|{\cal R}^0_0(\eps^{-2}n)|\geq 3 \big)=o(\eps)$, which is what condition (B2) in Theorem~\ref{T:convct1} amounts to.

\subsection{Proof of Lemma~\ref{L:ZBrownian}}

We follow the notation introduced in (\ref{truncation}). If paths in $\Gamma$ were Markovian in the sense that for any $u,v\in\Z$ with $u<v$, the law of $\Gamma_{u}^v$ depends only on $\Gamma_{u}(v)$ and not on the realization of paths in $\Gamma_{u}$ before time $v$, then Lemma~\ref{L:ZBrownian} would follow easily from Lemma~\ref{L:Zlocfin}, and condition (I) in Theorem~\ref{T:convct1} which was established in \cite[Prop.~3.3]{SS11}. Of course paths in $\Gamma$ are not Markovian. Nevertheless, Lemma~\ref{L:ZBrownian} asserts that the scaling limit of $\Gamma$ does satisfy the Markovian property described above, which is also a property satisfied by the Brownian web.
 To show how such Markovian property arises in the scaling limit, we will first remove
the dependence on the future by approximating $\Gamma_{u}(v)$ with ${\cal R}_{u}(v)$, where ${\cal R}$ is the set of right boundaries of all exploration clusters on $\Z^2_{\rm even}$. We will then approximate $\Gamma_{u}^v$ by $\Gamma_{v, {\cal R}_u}:=\{\gamma_z: z\in {\cal R}_{u}(v)\times\{v\}\}$, which removes the dependence on the past. Our proof of Lemma~\ref{L:ZBrownian} will consist of showing that the above approximations are accurate in the diffusive scaling limit.
\bigskip

\noindent
{\bf Proof of Lemma~\ref{L:ZBrownian}.} Given $t_0\in\R$, let $\eps_n\downarrow 0$ be such that, with ${\cal X}_n:= S_{\eps_n}\Gamma$, ${\cal X}_{n, t_0^-}$ converges weakly to ${\cal Z}$. First we show that it suffices to consider ${\cal X}_{n,0}$ in place of ${\cal X}_{n, t_0^-}$, which is just a technicality.

Using Skorohod's representation theorem, it is easily seen that
for any $\delta>0$, ${\cal X}_{n,t_0^-}^{t_0+\delta}$ converges weakly to ${\cal Z}^{t_0+\delta}$. Let $u_n=\lfloor \eps_n^{-1} t_0\rfloor$. Then ${\cal X}_{n, \eps_n u_n}$ is a.s.\ a closed subset of ${\cal X}_{n, t_0^-}$, and hence forms a tight sequence. Therefore by going to a subsequence if necessary, we can assume that $({\cal X}_{n, t_0^-}, {\cal X}_{n, \eps_n u_n})$ converges weakly to a limit $({\cal Z}, \tilde {\cal Z})$. Since ${\cal X}_{n,t_0^-}^{t_0+\delta}={\cal X}_{n, \eps_n u_n}^{t_0+\delta}$ because paths in $\Gamma$ coalesce when they intersect, we have $\tilde {\cal Z}^{t_0+\delta} = {\cal Z}^{t_0+\delta}$ a.s. In particular, to
identify the law of ${\cal Z}^{t_0+\delta}$ as coalescing Brownian motions starting from ${\cal Z}(t_0+\delta)\times\{t_0+\delta\}$,
we can replace the weakly convergent sequence ${\cal X}_{n,t_0^-}$ by ${\cal X}_{n, \eps_n u_n}$. Furthermore, because $u_n\in\Z$ and $\eps_n u_n\to t_0$, by the
translation invariance of $\Gamma$ under shifts by vectors in $\Z^2_{\rm even}$ and the a.s.\ equicontinuity of paths in $\tilde {\cal Z}$, we can take $t_0=0$ and
$u_n=0$ for all $n\in\N$. Lemma~\ref{L:ZBrownian} then reduces to showing that
\begin{itemize}
\item[($\dagger$)] For any $\delta>0$, if ${\cal X}_{n, 0}$ converges weakly to a limit ${\cal Z}$, then ${\cal Z}^\delta$ is distributed as $\Wi_{\delta, \cal Z}$, which is defined as in Lemma~\ref{L:ZBrownian}.
\end{itemize}
Note that ${\cal Z}(\delta)\times\{\delta\}$ is a.s.\ a locally finite subset of $\R^2$ by Lemma~\ref{L:Zlocfin}.

Let $v_n=\lfloor \eps_n^{-1}\delta\rfloor$. Note that the weak convergence of ${\cal X}_{n,0}$ to ${\cal Z}$ implies that
\be\label{Xn0delta}
({\cal X}_{n,0}, {\cal X}_{n,0}^{\eps_n v_n}) = S_{\eps_n}( \Gamma_0, \Gamma_0^{v_n}) \stackrel{\rm dist}{\Longrightarrow} ({\cal Z}, {\cal Z}^\delta) \qquad \mbox{as } n\to\infty.
\ee
Let ${\cal R}_0:=\{r_z: z=(x,0)\in\Z^2_{\rm even}\}$, and let $\Gamma_{v_n, {\cal R}_0} :=\{\gamma_z: z\in {\cal R}_0(v_n)\times\{v_n\}\}$. We will show that as $n\to\infty$,
\be\label{Lem4.2conv1}
S_{\eps_n} \Gamma_{v_n, {\cal R}_0}\ \stackrel{\rm dist}{\Longrightarrow}\ \Wi_{\delta, {\cal Z}},
\ee
and
\be\label{Lem4.2conv2}
S_{\eps_n}(\Gamma_0^{v_n}, \Gamma_{v_n, {\cal R}_0}) \stackrel{\rm dist}{\Longrightarrow} ({\cal Z}^\delta, {\cal Z}^\delta).
\ee
From (\ref{Lem4.2conv1}) and (\ref{Lem4.2conv2}), ($\dagger$) follows immediately.

The proof of (\ref{Lem4.2conv1}) is based on the fact that conditional on the realization of ${\cal R}_0(v_n)$, $\Gamma_{v_n, {\cal R}_0}$ is independent of what happens
before time $v_n$. To apply Skorohod's representation theorem later on, for each $n\in\N$, we construct $\Gamma^{[n]}_0$ and ${\cal R}^{[n]}_0$ from a random percolation edge configuration $\Omega^{[n]}$, such that $(\Gamma^{[n]}_0, {\cal R}^{[n]}_0)$ has the same distribution as that of $(\Gamma_0, {\cal R}_0)$. Let
$\gamma^{[n]}_z$, $r^{[n]}_z$, $\Gamma_0^{[n], v_n}$ and $\Gamma^{[n]}_{v_n, {\cal R}^{[n]}_0}$ denote the analogue of $\gamma_z$, $r_z$, $\Gamma_0^{v_n}$ and $\Gamma_{v_n, {\cal R}_0}$.

By \cite[Lemma~2.5]{SS11}, which controls the difference between $\gamma_z$ and $r_z$, there exists $C>0$ such that for each $\eps\in (0,1)$ and $N\in\N$,
$$
\begin{aligned}
\P\Big(\sup_{x\in [-\eps^{-1/2}N, \eps^{-1/2}N]\cap\Z} \ \sup_{t\in [0, 2\eps^{-1} \delta]} |r_{(2x,0)}(t) - \gamma_{(2x,0)}(t)| \geq \eps^{-1/4}\Big)
\leq 2\eps^{-1/2}N\ C \eps^{4}.
\end{aligned}
$$
Therefore by letting $N=\lceil\eps^{-1/2}\rceil$ and going to a subsequence of $(\eps_n)_{n\in\N}$, $\eps_n\downarrow 0$, if necessary, we can assume that
\be\label{intervalBC}
\sum_{n=1}^\infty \P\Big(\sup_{x\in [-\eps_n^{-1}, \eps_n^{-1}]\cap\Z} \ \sup_{t\in [0, 2\eps_n^{-1} \delta]} |r^{[n]}_{(2x,0)}(t) - \gamma^{[n]}_{(2x,0)}(t)| \geq \eps_n^{-1/4}\Big) < \infty.
\ee
By Borel-Cantelli, the events in the summation above occur a.s.\ only finitely many times regardless of how the sequence of percolation edge configurations $(\Omega^{[n]})_{n\in\N}$ are coupled.

By (\ref{Xn0delta}),
$$
S_{\eps_n}(\Gamma^{[n]}_0,\ \Gamma^{[n]}_{0}(v_n)\times\{v_n\}) \stackrel{\rm dist}{\Longrightarrow} ({\cal Z},\ {\cal Z}(\delta)\times\{\delta\})
\qquad \mbox{as } n\to\infty,
$$
where the second components are taken to be random variables taking values in the space of compact subsets of $\Rc$ (a suitable compactification of $\R^2$, see~\cite[Sec.~1.2]{SS11}), endowed with the Hausdorff topology. Using Skorohod's representation theorem to turn the above convergence into a.s.\ convergence, and then applying (\ref{intervalBC}) and Borel-Cantelli, we deduce that $S_{\eps_n}({\cal R}^{[n]}_0(v_n)\times\{v_n\})$ also converges weakly to ${\cal Z}(\delta) \times\{\delta\}$.

Let $\Omega^{[n]}_{(-\infty, v_n]}$ and $\Omega^{[n]}_{[v_n, +\infty)}$ denote respectively the configuration of edges in the percolation configuration $\Omega^{[n]}$ before and after time $v_n$.
Then ${\cal R}^{[n]}_0(v_n)$ depends only on $\Omega^{[n]}_{(-\infty, v_n]}$.  By Skorohod's representation theorem, we can couple $\Omega^{[n]}_{(-\infty, v_n]}$, $n\in\N$, such that a.s.\
\be\label{pointsetconv}
S_{\eps_n}({\cal R}^{[n]}_0(v_n)\times\{v_n\})\to {\cal Z}(\delta) \times\{\delta\}.
\ee
Since ${\cal Z}(\delta) \times\{\delta\}$ is a.s.\ locally finite by Lemma~\ref{L:Zlocfin}, we can label the points in ${\cal Z}(\delta) \times\{\delta\}$
successively by $(z_m)_{m\in\Z}$, where $z_m=(x_m,\delta)$ and $x_m<x_{m+1}$ for all $m\in\Z$. For each $n\in\N$ and $m\in\Z$, let $z_{n,m}^{\pm}:=(x_{n,m}^{\pm},v_n)$,
where
$$
\begin{aligned}
x_{n,m}^+ & :=  \max\big\{i\in {\cal R}_0(v_n): \frac{\sqrt{\eps_n}(i-\alpha v_n)}{\sigma} \leq \frac{x_m+x_{m+1}}{2} \big\}, \\
x_{n,m}^- & :=  \min\big\{i\in {\cal R}_0(v_n): \frac{\sqrt{\eps_n}(i-\alpha v_n)}{\sigma} \geq \frac{x_{m-1}+x_m}{2} \big\}.
\end{aligned}
$$
By (\ref{pointsetconv}), a.s.\ for each $m\in\Z$, $S_{\eps_n}z_{n,m}^{\pm}\to z_m$ as $n\to\infty$. Since $\Omega^{[n]}_{(-\infty, v_n]}$ and $\Omega^{[n]}_{[v_n, +\infty)}$ are independent, conditional on $\Omega^{[n]}_{(-\infty, v_n]}$, $n\in\N$, \cite[Prop.~3.3]{SS11} implies that
\be\label{zdpathconv}
S_{\eps_n}\big((\gamma_{z_{n,m}^+}, \gamma_{z_{n,m}^-})_{m\in\Z}\big)
\stackrel{\rm dist}{\Longrightarrow} \big((\Wi(z_m), \Wi(z_m))_{m\in\Z} \big),
\ee
where $(\Wi(z_m))_{m\in\Z}$ is a collection of coalescing Brownian motions starting at $(z_m)_{m\in\Z}$, embedded in a standard Brownian web $\Wi$ which is independent of $(z_m)_{m\in\Z}$.
Therefore applying Skorohod's representation theorem once more, conditional on
$\Omega^{[n]}_{(-\infty, v_n]}$, $n\in\N$, we can couple $\Omega^{[n]}_{[v_n,\infty)}$, $n\in\N$, such that the convergence in (\ref{zdpathconv}) becomes a.s.,
which implies that under such a full coupling of $\Omega^{[n]}$, $n\in\N$, we have
\be\label{zdpathconv2}
S_{\eps_n} \{\gamma_{z_{n,m}^{\pm}}: m\in\Z\} \to \{\Wi(z_m) : m\in\Z\}.
\ee
Our definition of $x_{n,m}^{\pm}$ guarantees that a.s.\ for each $m\in\Z$, if $n$ is sufficiently large, then
\be\label{R0vn}
\begin{aligned}
{\cal R}_0(v_n)\ \cap\ [x_{n,m}^-,\ x_{n,m}^+ ]\neq \emptyset \qquad \mbox{and} \qquad
{\cal R}_0(v_n)\ \cap\ (x_{n,m}^+,\ x_{n, m+1}^-)=\emptyset.
\end{aligned}
\ee
Since $\gamma_{(y, v_n)}$ is bounded between $\gamma_{z_{n,m}^-}$ and $\gamma_{z_{n,m}^+}$ for all $y\in [x_{n,m}^-, x_{n,m}^+]$,
(\ref{zdpathconv2}) and (\ref{R0vn}) imply that a.s.\
\be
S_{\eps_n} \Gamma_{v_n, {\cal R}_0} \to \{\Wi(z_m) : m\in\Z\}=\Wi_{\delta, \cal Z},
\ee
which proves (\ref{Lem4.2conv1}). Applying the same argument with $x_{n,m}^{\pm}$ replaced by $x_{n,m}^{\pm} \pm 2\lceil \eps_n^{-1/4}\rceil$ to accommodate
the difference between $\Gamma_0(v_n)$ and $\Ri_0(v_n)$, we deduce that $S_{\eps_n} \Gamma_0^{v_n} \to \Wi_{\delta, \cal Z}$ and thus
(\ref{Lem4.2conv2}), because points in $\Gamma_0(v_n)$ and ${\cal R}_0(v_n)$ a.s.\ approximate each other locally within a distance of $\eps_n^{-1/4}$ for all large $n$ by (\ref{intervalBC}) and Borel-Cantelli.
\qed

\section{Proof of Theorem~\ref{T:density}}\label{S:D}
The proof of Theorem~\ref{T:density} is essentially the same as that for \cite[Corollary~7.1]{NRS05}.
We will outline the main steps and point out the differences.

Analogous to \cite[Theorem~7.1]{NRS05}, by Theorem~\ref{T:main}, or rather, by the same proof as for Theorem~\ref{T:main}, we have
\be\label{Gamma0}
S_{\frac{1}{\sqrt n}} \Gamma_0 \Asto{n} \Wi_0,
\ee
where $\Gamma_0$ and $\Wi_0$ are defined from $\Gamma$ and $\Wi$ as in (\ref{truncation}), and $\Rightarrow$ denotes weak convergence of $\Hi$-valued random variables.
From (\ref{Gamma0}), it follows that
\be\label{Gamma0n}
(S_{\frac{1}{\sqrt n}} \Gamma_0)(1) \Asto{n} \Wi_0(1)
\ee
as random variables taking values in the space of compact subsets of $[-\infty, \infty]$, equipped with the Hausdorff topology. Since $\Wi_0(1)$ is a translation invariant point process on $\R$ with intensity $\frac{1}{\sqrt \pi}$, Fatou's lemma implies that
$$
\begin{aligned}
\frac{1}{\sqrt \pi}=\E\big[\,\big|\Wi_0(1)\cap [0,1]\big|\,\big] & \leq \liminf_{n\to\infty} \E\big[\,\big| (S_{\frac{1}{\sqrt n}} \Gamma_0)(1)\cap [0,1] \big|\,\big]\\
& = \liminf_{n\to\infty} \E[\, |\,\Gamma_0(n)\cap [\alpha n, \alpha n+\sigma\sqrt{n}]\,|\,] \\
& = \liminf_{n\to\infty} \frac{\sigma \sqrt{n}}{2} \P(\Gamma_0(n)\cap \{0,1\}\neq \emptyset),
\end{aligned}
$$
where we have used the translation invariance of $\Gamma_0(n)$. This gives the desired lower bound on $\P(\Gamma_0(n)\cap \{0,1\}\neq \emptyset)$.

The matching upper bound will follow from Lemma~\ref{L:densitydecay} and
\be\label{Ri0lim}
\lim_{n\to\infty} \frac{\sigma \sqrt{n}}{2} \P(\Ri_0(n)\cap \{0,1\}\neq \emptyset) = \frac{1}{\sqrt \pi},
\ee
where $\Ri_0$ is defined as in Lemma~\ref{L:densitydecay}. To prove (\ref{Ri0lim}), we observe that (\ref{Gamma0}) and (\ref{Gamma0n}) also hold with $\Gamma_0$ replaced by $\Ri_0$, because $\Ri_0$ has non-crossing paths by \cite[Lemma~3.1]{SS11}, condition (I) in Theorem~\ref{T:convct1} was verified in~\cite[Prop.~3.3]{SS11} for $S_\eps\Ri$ as well as $S_\eps\Gamma$, and conditions (B1) and (B2) in Theorem~\ref{T:convct1} were verified in~\cite{SS11} for $S_\eps\Gamma$ by first replacing it with
$S_\eps\Ri$ (alternatively, note that the verification of condition (E) in Section~\ref{S:E} applies to $S_\eps\Ri$ as well).

We can then proceed as in~\cite[Theorem~7.3]{NRS05} to strengthen the weak convergence of $(S_{1/\sqrt n}\Ri_0)(1)\Rightarrow\Wi_0(1)$ by regarding $(S_{1/\sqrt n} \Ri_0)(1)$ and $\Wi_0(1)$ as random counting measures on $\R$ with convergence with respect to the vague topology, where each point in $(S_{1/\sqrt n} \Ri_0)(1)$ and $\Wi_0(1)$ is replaced by a delta measure at that point. Note that $(S_{1/\sqrt n} \Ri_0)(1)$, $n\in\N$, is a tight family of random counting measures by the crude density bound in Lemma~\ref{L:densitydecay}. Furthermore, any subsequential weak limit of $(S_{1/\sqrt n} \Ri_0)(1)$, $n\in\N$, is necessarily a simple point process. The proof of this fact follows the same steps as in the proof of~\cite[Theorem~7.3]{NRS05}, once we observe that $\Ri_0(n)$ satisfies the negative correlation inequality
\be\label{negcor}
\P(i, j\in \Ri_0(n)) \leq \P(i \in \Ri_0(n)) \P(j \in \Ri_0(n)) \quad \forall\,  i<j \mbox{ with } (i,n), (j,n)\in\Z^2_{\rm even}.
\ee
Indeed, when $\{i\in \Ri_0(n)\}$ and $\{j\in \Ri_0(n)\}$ both occur, there must exist $z_1, z_2\in \Z\times\{0\}$ such that $r_{z_1}(n)=i$ and $r_{z_2}(n)=j$. In particular, the two associated exploration clusters $C_{z_1}(n)$ and $C_{z_2}(n)$ must be disjoint, and hence $\{i\in \Ri_0(n)\}$ and $\{j\in \Ri_0(n)\}$ occur disjointly, i.e., occur using two disjoint sets of oriented edges. Therefore by Reimer's inequality~\cite{R00}, these two events are negatively correlated, which gives
(\ref{negcor}). Together with the weak convergence of $S_{1/\sqrt n}\Ri_0\Rightarrow\Wi_0$, we can then deduce the weak convergence of $(S_{1/\sqrt n}\Ri_0)(1)\Rightarrow\Wi_0(1)$ as random counting measures, following the same arguments as in the proof of~\cite[Theorem~7.3]{NRS05}. Applying the same argument as
in the proof of~\cite[Corollary~7.1]{NRS05} then gives (\ref{Ri0lim}), which concludes the proof of Theorem~\ref{T:density}.
\qed
\smallskip

\noindent
{\bf Remark.} Note that in the proof of Theorem~\ref{T:density}, we cannot apply Reimer's inequality to deduce directly the analogue of (\ref{negcor}) for $\Gamma_0(n)$ because of the dependence of $\Gamma_0(n)$ on the future, which is why we switched from $\Gamma_0(n)$ to $\Ri_0(n)$ instead.
\bigskip

\noindent
{\bf Acknowledgement} A. Sarkar wishes to thank NWO and Vrije Universiteit Amsterdam for support. R.~Sun is supported by grant R-146-000-119-133 from the National University of Singapore.

\end{document}